%% file: 0.tex
\documentclass[10pt]{article}

\usepackage{amsmath}
\usepackage{amsthm}
\usepackage{euscript}

\usepackage{diagrams}

\usepackage[all]{xy}

\input macros.tex
\input theorems.tex

\title{Logical Construction of Final Coalgebras\footnote{Research
    supported by the European Community through a Marie Curie
    Individual Fellowship.}}

\author{%
  Luigi Santocanale\\%
  Laboratoire d'Informatique Fondamentale\\%
  Centre de Math\'ematiques et Informatique \\%
  39, rue Joliot-Curie - F-13453 Marseille Cedex13 \\%
  {\tt luigi.santocanale@cmi.univ-mrs.fr}\\%
}
\date{\today}

\begin{document}

\maketitle

\input abstract.tex


\input intro.tex
\input preliminaries.tex
\input automata.tex
\input characterization.tex
\input specialQaut.tex 
\input further.tex

\input ack.tex


\input 0.bbl\end{document}

%% file: macros.tex
\newcommand{\hspp}{\hspace{1mm}}

\newcommand{\bbt}{\mathsf{T}}
\newcommand{\N}{\mathsf{N}}
\newcommand{\die}[1]{#1{}^{\#}}
\newcommand{\be}[1]{#1{}^{\flat}}
\newcommand{\id}{\mathtt{id}}
\newcommand{\comp}{\cdot}
\newcommand{\ev}{\mathtt{ev}}
\newcommand{\iso}{\cong}
\newcommand{\pos}[1][\hspace{0.5em}]{\langle#1\rangle}

\newcommand{\inj}{{\mathtt{in}}}
\newcommand{\set}[1]{\{ #1\}}
\newcommand{\Cat}[1]{{\EuScript{#1}}}
\newcommand{\pb}[2][0]{
  \begin{array}[t]{@{\hspace{0mm}}c@{\hspace{0mm}}}
    \times \\[-#1mm]
    \scriptstyle{#2}
  \end{array}
}

\newcommand{\F}{\widehat{F}}
\newcommand{\s}{\widehat{s}}
\newcommand{\G}{\widetilde{G}}

\newcommand{\fold}[1]{\{\!|#1|\!\}}

\newcommand{\zero}{\mathtt{z}}
\renewcommand{\succ}{\mathtt{s}}
\newcommand{\mul}{\mathtt{m}}
\newcommand{\inF}{\mathtt{i}}

\newcommand{\head}{\mathtt{h}}
\newcommand{\tail}{\mathtt{t}}
\newcommand{\comul}{\mathtt{d}}
\newcommand{\prG}{\mathtt{p}}

\newcommand{\Pautomata}{\mathrm{Aut}(P)}
\newcommand{\Qautomata}{\mathrm{Aut}(Q)}
\newcommand{\dQautomata}{\mbox{$\delta$-$\mathrm{Aut}(Q)$}}
\newcommand{\coalg}[1]{\mathrm{CoAlg}(#1)}
\newcommand{\Pcoalg}{\mathrm{CoAlg}(P)}
\newcommand{\Qcoalg}{\mathrm{CoAlg}(Q)}
\newcommand{\Falg}[1][]{\mathrm{Alg}(F)_{#1}}
\newcommand{\Gcoalg}{\mathrm{CoAlg}(G)}
\newcommand{\inje}[2][]{\iota^{#2}_{#1}}

\newcommand{\Aut}[1]{\EuScript{#1}}

\newcommand{\putequation}[1]{\refstepcounter{equation}{\protect\label{#1}}\mathrm{(\theequation)}}
\newcommand{\by}[1]{\tag*{\textrm{by (\ref{#1})}}}

\newgraphescape{E}[1]{[]!{c+(#1,-0.5)}}

\newcommand{\mydiagram}[2][6em]{
  \mbox{$\begin{array}{@{\hspace{0mm}}c@{\hspace{0mm}}}
      \xy\xygraph{!{<#1,0cm>:<0cm,#1>::}#2}\endxy
    \end{array}$}
  }

\newgraphescape{W}{[](
  [d(0.10)r(0.25)]-[d(0.15)]="MyspecialEnd"
  [r(-0.15)]-"MyspecialEnd"
  )
}

\newgraphescape{Z}{[](
  [d(0.10)r(0.25)]-[d(0.10)][d(0.05)]
  [l(0.05)]="MyspecialEnd"
  [r(-0.10)]-"MyspecialEnd"
  )
}

\newgraphescape{s}[6]{
  []*+{#1}="1"
  [d(#5)]*+{#2}="2"
  [r(#6)]*+{#4}="4"
  "1"[r(#6)]*+{#3}="3"
}

\newgraphescape{a}[4]{
  "1":"2"#1
  "1":"3"#2
  "2":"4"#3
  "3":"4"#4
}



%% file: theorems.tex
\theoremstyle{plain}
\newtheorem{thrm}{Theorem}[section]
\newtheorem{lmm}[thrm]{Lemma}
\newtheorem{crllr}[thrm]{Corollary}
\newtheorem{prpstn}[thrm]{Proposition}

\theoremstyle{definition}
\newtheorem{dfntn}[thrm]{Definition}

\newtheorem{xmpl}[thrm]{Example}

\theoremstyle{plain}

\newenvironment{theorem}{\begin{thrm}}{\end{thrm}}
\newenvironment{proposition}{\begin{prpstn}}{\end{prpstn}}
\newenvironment{lemma}{\begin{lmm}}{\end{lmm}}
\newenvironment{corollary}{\begin{crllr}}{\end{crllr}}
\newenvironment{example}{\begin{xmpl}}{\end{xmpl}}
\newenvironment{definition}{\begin{dfntn}}{\end{dfntn}}



%% file: abstract.tex
\begin{abstract}
We prove that every finitary polynomial endofunctor of a category
$\Cat{C}$ has a final coalgebra if $\Cat{C}$ is locally Cartesian
closed, has finite disjoint coproducts and a natural number object.
More generally, we prove that the category of coalgebras for such an
endofunctor has all finite limits.
\end{abstract}


%% file: intro.tex
\section*{Introduction}

A goal of this paper is to prove that every polynomial endofunctor
\begin{align}
  \label{def:P}
  P(X) & = \sum_{i = 1, \ldots ,n} \Omega_{i} \times X^{A_{i}}
\end{align}
of a given category $\Cat{C}$ has a final coalgebra, assuming that (1)
$\Cat{C}$ is locally Cartesian closed, (2) it has finite disjoint
coproducts, (3) it has a natural number object. This statement follows
from a more general result which also shows that the category of
$P$-coalgebras has all finite limits. An immediate consequence of this
statement is that the functor $P$ is completely iterative in the sense
of \cite{aczeletal} and generates a cofree comonad.  To this end
observe that the collection of polynomial endofunctors is closed
w.r.t. addition and multiplication by a fixed object $A$ of
$\Cat{C}$ (the latter up to natural isomorphism), and recall that a
comonad cofreely generated by $P$ is essentially the same as the
parametrized final coalgebra of $A \times P(X)$.

\hspp

Our proof is inspired by a particular set theoretic representation --
see \cite{courcelle,esikbloom} for example -- of the final
$P$-coalgebra
which we illustrate next.
\begin{example}
  \label{ex:main}
  Consider the functor $P(X) = \Omega_{1} \times X + \Omega_{2}\times
  X^{2}$, where $\Omega_{1} = \set{f}$ and $\Omega_{2} = \set{g}$.
  Its final coalgebra is the set of finite and infinite terms over the
  signature $\set{f,g} = \Omega_{1} + \Omega_{1}$, $f$ a unary
  function symbol and $g$ a binary function symbol.  Let now $A =
  \set{(f,1),(g,1),(g,2)}$ be the disjoint sum of the arities of the
  signature and let $\Omega = \set{f,g ,\bot}$ be obtained from the
  previous signature by addition of the new symbol $\bot$.  The
  infinite terms over the signature $\set{f,g}$ are in bijection with
  those complete $A$-branching trees labeled in $\Omega$ satisfying
  the following clauses:
\begin{enumerate}
\item \label{clause:f} If a node is labeled by $f$, then the son in
  the direction $(f,1)$ is not labeled by the symbol $\bot$ while all
  the sons in the directions $(g,1),(g,2)$ are labeled
  by $\bot$.
\item \label{clause:g} If a node is labeled by $g$, then the son in
  the direction $(f,1)$ is labeled by the symbol $\bot$ and the sons
  in the directions $(g,1),(g,2)$ are not labeled by $\bot$.
\item \label{clause:notbot} The root of the tree is not labeled by the
  symbol $\bot$.
\end{enumerate}
A complete tree satisfying these conditions is represented in figure
\ref{fig:tree}.
\begin{figure}[h]
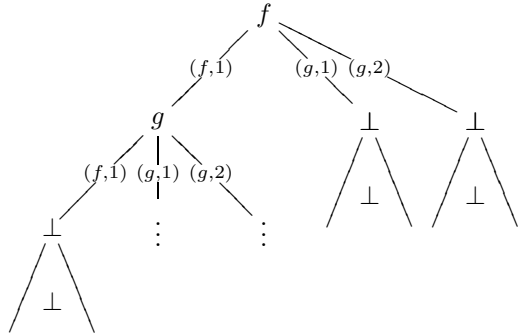

  \newgraphescape{v}{ []*+{\bot}
    (-[l(0.4)d(1)],[d(0.7)]*+{\bot},-[r(0.4)d(1)]) }
  $$
  \mydiagram[4em]{ []*+{f} ( -[ld]*+{g}|{(f,1)}
    (-[ld]*+{}|{(f,1)}!v, -[d]*+{\vdots}|{(g,1)},
    -[rd]*+{\vdots}|{(g,2)} ), -[rd]*+{\bot}|{(g,1)}!v,
    -[rrd]*+{\bot}|{(g,2)}!v, ) }
  $$
  \caption{Infinite terms as complete trees}
  \label{fig:tree}
\end{figure}
\end{example}
This example can obviously be generalized to arbitrary set-theoretic
polynomial endofunctors and arbitrary signatures. If we let $A =
\sum_{i} A_{i}$ and $\Omega = \set{\bot} + \sum_{i} \Omega_{i}$, the
set-theoretic final $P$-coalgebra is a subset of the function space
$\Omega^{A^{\ast}}$.  This function space can also be seen as the set
of complete $A$-branching trees with labels in $\Omega$; as such it
carries the structure of a final coalgebra for the functor $Q(X) =
\Omega \times X^{A}$.

\hspp

We want to investigate the process by which a final $P$-coalgebra is
extracted from a final $Q$-coalgebra by means of logical operations.
The outcome of our investigation can be synthesized as follows.  A
full and faithful `completion' functor $K$ from the category of
$P$-coalgebras to the category of $Q$-coalgebras is defined.  This
functor factors as a full and faithful right adjoint $K^{+}$ followed
by a full and faithful left adjoint $I$.  This suffices to argue that
the category of $P$-coalgebras has finite limits provided that the
category of $Q$-coalgebras has finite limits -- a standard lemma in
the theory of reflective and coreflective subcategories.  Also, under
our assumptions, the category of $Q$-coalgebras has indeed all finite
limits. 

This extraction process can be carried out by means of the weak
internal logic corresponding to the categorical properties (1)-(3).
This is the minimal logic in which the clauses of Example
\ref{ex:main} are expressible. For example, the use there of classical
logic is restricted to the constructive Boolean logic of extensive
categories, see \cite{cockettlack}.  Interestingly, clauses
(\ref{clause:f})-(\ref{clause:g}) point to the existence of the left
adjoint $I$ while clause (\ref{clause:notbot}) of \ref{ex:main} is
related to the right adjoint $K^{+}$.

\hspp

Our interest in this problem is part of a general investigation
\cite{indcoind,paritygames} on the relationships between induction and
coinduction, i.e.  initial algebras and final coalgebras of functors.
In \cite{indcoind} we proved that, given an adjoint pair of
endofunctors $F \dashv G$, if a free $F$-algebra functor $\F$ exists
and has a right adjoint $\G$, then $\G$ is a cofree $G$-coalgebra
functor.\footnote{The converse holds too: if both a free $F$-algebra
  functor $\F$ and a cofree $G$-coalgebra functor $\G$ exist, then
  $\F$ is left adjoint to $\G$, see
  \cite{arbibmanes,altenkirch,indcoind}.}  For example, if $A^{\ast}$
is the free monoid generated by $A$ in a Cartesian closed category,
then $A^{\ast} \times X$ is the free $A \times (-)$-algebra generated
by $X$, and the function space $X^{A^{\ast}}$ turns out to be the
cofree $(-)^{A}$-coalgebra over $X$. Therefore if a locally Cartesian
closed category has a natural number object $\N$, then it has the free
monoids $A^{\ast}$ -- they are computed as partial
products\footnote{See \cite{dyckofftholen}.} of $A$ with the arrow $\N
\times \N \rTo^{+} \N\rTo^{s} \N$ -- and consequently it has final
$(-)^{A}$-coalgebras.  Briefly, some final coalgebras arise from
duality.  

A question suggested from the set theoretic example but left
open in \cite{indcoind} was whether and under which conditions the
duality is the starting point for constructing other final coalgebras.
The following result, expressed within the framework of
\cite{indcoind}, gives a positive answer to this question: given
adjoint pairs $F_{i} \dashv G_{i}$, $i = 1,\ldots ,n$, the category of
coalgebras for the endofunctor
\begin{align*}
  P(X) & = 
  \sum_{i = 1, \ldots ,n} \Omega_{i} \times G_{i}(X)
\end{align*}
has all finite limits under the assumption that $\Cat{C}$ has the
properties (1)-(2) and that both a free algebra functor of $\sum_{i} F_{i}$ and a cofree coalgebra functor of $\prod_{i} G_{i}$ exist.  We recover the previous statement
letting $F_{i} X = A_{i} \times X$, $G_{i} X = X^{A_{i}}$.
Using this framework it is possible to argue that functorial
systems of equations such as
$$
\left
  \{ X_{k} 
  \;=\;
  \sum_{i= 1,\ldots ,n} \Omega_{i,k} \times X_{f(k)}^{A_{i,k}}
\right \}_{k = 1,\ldots ,m}
$$
have a greatest solution. To this end observe that conditions
(1)-(2) are stable under formation of products of categories and use
results from \cite{cockettjay}.

\hspp

A motivation for this work has been that of deriving existence of
final coalgebras on axiomatic bases, without either imposing strong
completeness properties on our models, nor relying on the full topos
theoretic categorical structure. As a matter of fact the construction
presented here depends only on a restricted fragment of this
structure: we do not need a subobject classifier nor a factorization
system corresponding to existential quantification.  

In concrete categories -- mainly locally presentable categories
\cite{lpac} -- final coalgebras are usually constructed as limits of
the sequence of iterated functorial applications beginning at the
terminal object \cite{barr,adgr,worrel}.  For these categories final
coalgebras of polynomial functors can be constructed in this way.  On
the other hand, several categories studied in computer science lack
the completeness properties required to successfully perform the
terminal sequence construction.  A priori, we list among them the
effective topos \cite{hyland} and the free topos generated by a
countable language \cite{lambekscott}. Yet, elementary toposes admit
final coalgebras of partial product functors, thus of polynomial
functors: in \cite{johnstoneworrel} their final coalgebras are
constructed as internal limits of the terminal sequence. Existence of
these internal limits depends on the development of internal category
theory and on the theory of iterative data in an elementary topos
\cite{johnstonewraith}. 

The full topos theoretic structure is often considered too strong, in
particular when designing programming languages with a categorical
semantics in mind -- see \cite{charityI,charityII} and the language
Charity for an example -- or considering categorical universes of
predicative mathematics \cite{moerdijkpalmgren}.  Locally Cartesian
closed categories with disjoint coproducts have often been
proposed as an alternative choice \cite{cockettjay,containers}.  
The reason for not assuming completeness stems from considerations on
initial semantics of typed programming languages. These are usually
interpreted in categories with some kind of structure, and, among
these structured categories, the initial one plays the role of a
canonical model. This model will inevitably lack the completeness
properties: for example, a natural number object in such a category
cannot be a countable coproduct of the terminal object, since the
homsets are at most countable.\footnote{This remark is due to R.
  Cockett.}  We finally mention that recursive-theoretic categories
such as $\omega$-$Sets$ and $PER$ are locally Cartesian closed without
being complete or topoi \cite[\S 1.2.7]{jacobs}.


\hspace{1.5mm}

The paper is structured as follows: in section \ref{sec:prel} we
overview the mathematical setting and introduce the notation. In
section \ref{sec:automata} we represent the categories of
$P$-coalgebras and $Q$-coalgebras in an `automatic' form, exploiting
extensiveness of $\Cat{C}$. This allows to easily define the
completion functor $K$, its factorization, and to argue that its first
factor $K^{+}$ has a left adjoint.  In section
\ref{sec:characterization} we make use of locally Cartesian closeness
to internally characterize the domain of the second factor $I$. In the
sequent section \ref{sec:coreflection} we show that the factor $I$ has
a right adjoint. We add concluding remarks in section
\ref{sec:further}.


%% file: preliminaries.tex
\section{Preliminaries}
\label{sec:prel}

The goal of this paper is to prove the following statement:
\begin{theorem}
  \label{main:theorem}
  Let $\Cat{C}$ be a locally Cartesian closed category with finite
  disjoint coproducts.  Let $F_{i} \dashv G_{i}$, $i \in I$, be a
  finite collection of adjoint endofunctors of $\Cat{C}$ and define
  \begin{align*}
    F(X) & = \sum_{i \in I} F_{i}(X)\,,
    &
    G(X) & = \prod_{i \in I} G_{i}(X)\,.
  \end{align*}
  If a pair of adjoint endofunctors $\F\dashv\G$ -- where $\F$ is a
  free $F$-algebra functor or equivalently $\G$ is a cofree
  $G$-coalgebra functor -- exists, then the functor
  \begin{align}
    \label{eq:defP}
    P(X) & = 
    \sum_{i \in I} \Omega_{i} \times G_{i}(X)
  \end{align}
  has a final coalgebra. More generally, the category of
  $P$-coalgebras has all finite limits.
\end{theorem}
By examining the proof the reader will convince himself that if
$\Cat{C}$ has enough completeness properties, then the statement holds
for an indexing set $I$ of a given infinite cardinality.  We begin
explaining the statement and recalling some basic results we will need
later.

A category $\Cat{C}$ is \emph{locally Cartesian closed} if it has a
terminal object $1$ and each slice category $\Cat{C}/C$ is a Cartesian
closed category \cite{day,penon,seely}.  Hence $\Cat{C}$ is itself
Cartesian closed, being equivalent to the slice $\Cat{C}/1$.  If
restricted to monic arrows with same codomain $C$, the local Cartesian
closed structure endows the collection of all subobjects of $C$ with
the structure of a Brouwerian semilattice: it has finite intersections
$\bigwedge$ and an implication operation $\rightarrow$ satisfying usual
axioms. Moreover pulling back along an arrow preserves this structure.

A locally Cartesian closed $\Cat{C}$ is a \emph{distributive}
category, meaning that if we construct the pullbacks
$$
\mydiagram{
 []!s{P_{i}}{B_{i}}{A}{\sum_{i \in I} B_{i}}{1}{1.5}
 !a{}{^{\pi_{i}}}{^{\inj_{i}}}{^{f}}
 "1"!W
}
$$
where the $\inj_{i}: B_{i} \rTo \sum_{i} B_{i}$ are coproduct
injections, then the diagram $(\pi_{i}: P_{i} \rTo A)_{i \in I}$ is
again a coproduct. $\Cat{C}$ is an \emph{extensive} category if the
converse condition holds: if $(\pi_{i}: P_{i} \rTo A)_{i \in I}$ is a
coproduct diagram and the diagram above commutes, then it is a
pullback.  Rephrased, coproduct injections are Cartesian natural
transformations.  A locally Cartesian closed category has finite
disjoint coproducts (coproducts are disjoint if the intersection of
distinct coproduct injections is an initial object) if and only it has
finite coproducts and is an extensive category.\footnote{In the
  following we will only use the fact that coproduct injections are
  Cartesian natural transformation.}  The following property of
distributive categories, see \cite{robin}, will be frequently used:
\begin{lemma}
  \label{lemma:robin}
  In a distrubutive category coproduct injections are monic.
\end{lemma}

An \emph{$F$-algebra} is a pair $(Q,s)$ with $s : F Q \rTo Q$. A
morphism of $F$-algebras from $(Q,s)$ to $(Q',s')$ is an arrow $f : Q
\rTo Q'$ such that $s \comp f = Ff \comp s'$.  $F$-algebras and their
morphisms form a category $\Falg$. The category $\Gcoalg$ of
\emph{$G$-coalgebras} is defined dually. By saying that a \emph{free
  $F$-algebra functor} exists we mean that the forgetful functor
$U_{F} : \Falg \rTo \Cat{C}$, sending $(Q,s)$ to $Q$, has a left
adjoint.  Spelled out, for every object $X$ we can find an object $\F
X$ and a diagram
\begin{align*}
  \zero_{X} &: X \rTo \F X 
  & 
  \succ_{X} & : F\F X \rTo \F X
\end{align*}
with the initial property w.r.t. similar diagrams: for every pair
$(a,f)$, where $a: X \rTo A$ and $f : FA \rTo A$, there exists a
unique arrow $\fold{a,f}: \F X \rTo A$ such that
\begin{align*}
  \zero_{X} \comp \fold{a,f} & = a\,,
  & \succ_{X} \comp \fold{a,f} & = F
  \fold{a,f} \comp f \,.
\end{align*}
Similarly, by saying that a cofree $G$-coalgebra functor exists, we
mean that the forgetful functor $U_{G} : \Gcoalg \rTo \Cat{C}$ has a
right adjoint. Spelled out, for every object $X$ we can find an object
$\G X$ and a diagram
\begin{align*}
  \head_{X} &: \G X \rTo  X 
  & 
  \tail_{X} & : \G X \rTo G \G X
\end{align*}
with the  final property w.r.t. similar diagrams. Clearly, $\F$ and
$\G$ are functors, obtained by composing the left adjoint with $U_{F}$
and the right adjoint with $U_{G}$; $\zero,\succ, \head, \tail$ are
natural transformations.

If a free $F$-algebra functor $\F$ is given and we define
\begin{align*}
  \inF_{X} & = F\zero_{X} \comp \succ_{X} : FX \rTo \F X\,,
  &
  \mul_{X} & = \fold{\id_{\F X},\succ_{X}} : \F \F X \rTo \F X\,,
\end{align*}
then the tuple $\langle \F, \inF, \zero, \mul\rangle$ is the
\emph{free monad} generated by $F$ \cite{barrfreetriples,adfreemonad}.
Dually, if a cofree $G$-coalgebra functor $\G$ is given, then we can
define $\prG_{X} : \G X \rTo G X$ and $\comul_{X} : \G X \rTo \G \G X
$ so that $\langle \G, \prG,\head,\comul\rangle$ is the cofree comonad
generated by $G$.  Let $F \dashv G$ be a pair of adjoint endofunctors
of an arbitrary category. In \cite{indcoind} we proved the following
facts. If $\F$ is a free $F$-algebra functor and $\G$ is a cofree
$G$-coalgebra functor, then $\F$ is left adjoint to $\G$ -- see also
\cite{arbibmanes,altenkirch}.  Conversely, if a free $F$-algebra
functor $\F$ is given and has a right adjoint $\G$, then $\G$ can be
endowed with the structure of a cofree $G$-coalgebra functor. Dually:
if a cofree $G$-coalgebra functor $\G$ is given and has a left adjoint
$\F$, then $\F$ can be can be endowed with the structure of a free
$F$-algebra functor. The proof of these fact relied on the following
well known isomorphisms, see \cite{arbibmanes,barrfreetriples} and
\cite[\S 3.7]{barrwells}:
\begin{lemma}
  \label{lemma:isos}
  The category $\Falg$ of algebras for an endofunctor $F$ is isomorphic
  to:
  \begin{enumerate}
  \item \label{item:algcolag} the category $\coalg{G}$ for a functor
    $G$ right adjoint to $F$;
  \item the Eilenberg-Moore category of algebras for
    $(\F,\zero,\mul)$, the free monad generated by $F$;
  \item if $\G$ is right adjoint to $\F$, the Eilenberg-Moore category
    of coalgebras for the comonad on $\G$ dual to $(\F,\zero,\mul)$.
  \end{enumerate}
  These isomorphisms commute with the respective forgetful functors.
\end{lemma}
\noindent
We recall the nature of these isomorphisms. The first isomorphism
sends an algebra $s:FX \rTo X$ to its transpose $\be{s}: X \rTo GX$.%
\footnote{We shall use in this paper the notation $\be{(\cdot)}$ for
  the transpose of arrows under the adjunction, with $\die{(\cdot)}$
  for the inverse correspondence.}  The second isomorphism is obtained
by sending an algebra $s$ to $\s = \fold{\id_{X},s}$, the unique
arrow $\s: \F X \rTo X$ such that $\zero_{X} \comp \s = \id_{X}$ and
$\succ_{X} \comp \s = F \s \comp s$; it is easily seen that this is an
algebra for the  monad on $\F$.  Finally, if $(\F,\zero,\mul)$ is
a monad and $\G$ is right adjoint to $\F$ with $\ev$ the counit of the
adjunction, the dual comonad on $\G$ is defined as follows: $\head_{X}
= \zero_{\G X} \comp \ev_{X}$ is the counit of the cofree comonad
while the comultiplication $ \comul_{X}: \G X \rTo \G \G X $ is
obtained by transposing twice the arrow
$$
\mul_{\G X} \comp \ev_{X} : \F\F \G X \rTo \F \G X \rTo X\,.
$$
The isomorphism between algebras and coalgebras of
\ref{lemma:isos}.\textit{i} restricts to an isomorphism
between Eilenberg-Moore algebras and Eilenberg-Moore coalgebras.

We shall use the following fact: given an endofunctor $F$ of a
category $\Cat{C}$, we can define a new functor from the slice
category $\Cat{C}/C$ to the slice category $\Cat{C}/FC$: an object
$(X,x)$ -- where $x : X \rTo C$ -- is sent to $(FX,Fx)$ and an arrow
$f : X \rTo Y$ such that $f\comp y = x$ is sent to $Ff$. The following
lemma is well known \cite[\S 4.1.3]{elephant}.
\begin{lemma}
  \label{lemma:univquantification}
  If $F$ has a right adjoint $G$ and $\Cat{C}$ has pullbacks then the
  functor $F: \Cat{C}/C \rTo \Cat{C}/FC$ has also a right adjoint
  $\forall_{F}$ which is computed by pulling back along the unit of
  the adjunction:
  $$
  \mydiagram[6em]{
    [](!W)!s{\forall_{F} Q }
    {C}{GQ}{GFC}
    {1}{1.5}
    !a{^{\forall_{F}q}}{}{^{\eta_{C}}}{^{Gq}}
  }
  $$
\end{lemma}
\noindent
Observe that if $q$ is monic then $Gq$ is also monic, and therefore
$\forall_{F}q$ is monic.

We conclude this section introducing the notation. In the statement of
Theorem \ref{main:theorem} the functor $F$ is the sum $ \sum_{i \in I}
F_{i}$; hence we shall consider several restrictions of the
multiplication $\mul$ and of the action $s: F X \rTo X$ of an
arbitrary $F$-algebra:
\begin{align*}
  \mul_{i} & = \inj_{F_{i}\F} \comp \inF_{\F} \comp \mul 
  & s_{i} & = \inj_{F_{i}X} \comp \inF_{X} \comp \s \\
  & =
  \inj_{F_{i}\F} \comp \succ_{\F} : F_{i}\F  \rTo \F \,,
  &
  & =
  \inj_{F_{i}X} \comp s_{X} : F_{i}X  \rTo X \,.
\end{align*}
We shall occasionally say that $X$ is a subobject of $Y$, in which
case we shall use $\inje[Y]{X}$ for the intended monic arrow $X \rTo
Y$; we shall only write $\inje{X}$ if $Y$ is understood.


%% file: automata.tex
\section{Coalgebras as Automata, a Reflector}
\label{sec:automata} 
 
The category $\Cat{C}$ being extensive, we can represent coalgebras of
polynomial endofunctors by means of tuples, in an analogous way we
usually represent coalgebras as automata in the category of sets and
functions.  We make explicit this correspondence in this section: we
define a category $\Pautomata$ of $P$-automata and argue it is
equivalent to the category of $P$-coalgebras, $P$ being the functor
defined in (\ref{eq:defP}).  Similarly we define a category
$\Qautomata$ of $Q$-automata for the functor
\begin{align*}
  Q(X) & =
  \Omega \times \prod_{i \in I} G_{i}(X)\,, 
\end{align*}
where $\Omega$ is $1 + \sum_{i \in I} \Omega_{i}$, and we argue it is
equivalent to the category of $Q$-coalgebras.  We define then a
`completion' functor
$$
K : 
 \Pautomata
\rTo \Qautomata 
$$
and begin its study. We shall eventually see that this functor is
responsible for the existence of a terminal object and finite limits
in $\Pautomata$ -- and therefore in $\Pcoalg$.

We introduce here the following notation: we let $J = \set{0} \cup I$,
with $0 \not\in I$, $\Omega_{0} = 1$ is a terminal object, and for $K
\subseteq J$ we let $\Omega_{K} = \sum_{k \in K} \Omega_{k}$, so that
$\Omega = \Omega_{J} = 1 + \Omega_{I}$. A similar notation is used for
arbitrary families $\set{Q_{j}}_{j \in J}$, for example $Q_{I}$ stands
for $\sum_{i \in I} Q_{i}$.

\begin{definition}
  The category $\Pautomata$ of \emph{$P$-automata} is defined as
  follows:
  \begin{itemize}
  \item An object of $\Pautomata$ is a tuple 
    \begin{align*}
      \Aut{A} & 
      = \langle \,
      \{\inj_{i}:Q_{i}\rTo Q\}_{i \in I},
      \{s_{i}\}_{i \in I}, \{h_{i}\}_{i \in I}\,
      \rangle
    \end{align*}
  where:
  \begin{itemize}
  \item $\{\inj_{i}:Q_{i}\rTo Q\}_{i \in I}$ is a coproduct diagram
      in $\Cat{C}$,
    \item  $s_{i}: F_{i}Q_{i} \rTo Q$ for each $i \in I$,
    \item $h_{i} : A_{i} \rTo \Omega_{i}$ for each $i \in I$.
    \end{itemize}
  \item An arrow $f$ from $\Aut{A}$ to $\Aut{B}$ is an arrow $f :
    Q^{\Aut{A}} \rTo Q^{\Aut{B}}$ in $\Cat{C}$ such that:
    \begin{itemize}
    \item $\inj^{\Aut{A}}_{i}\comp f$ factors (uniquely because of
      Lemma \ref{lemma:robin}) through $\inj^{\Aut{B}}_{i}$:
      \begin{align*}
        \inj^{\Aut{A}}_{i}\comp f
        & = f_{i} \comp \inj^{\Aut{B}}_{i}\,,
      \end{align*}
    \item for each $i \in I$ the following equations hold:
      \begin{align*}
        s^{\Aut{A}}_{i}\comp f & = F_{i}f_{i}
        \comp s^{\Aut{B}}_{i}\,, &
        h^{\Aut{A}}_{i} & = f_{i} \comp h^{\Aut{B}}_{i} \,.
      \end{align*}
    \end{itemize}
  \end{itemize} 
\end{definition}
We shall occasionally use the standard terminology for automata: $Q$
is the carrier of $\Aut{A}$, the $s_{i}$ are actions, and the $h_{i}$
are labeling arrows.
\begin{definition}
  The category $\Qautomata$ of \emph{$Q$-automata} is defined as
  follows:
  \begin{itemize}
  \item An object of $\Qautomata$ is a tuple 
    \begin{align*}
      \Aut{A} & 
      = \langle \,
      \{\inj_{j}:Q_{j}\rTo Q\}_{j \in J},
      \{s_{i,j}\}_{i \in I, j \in J},
      \{h_{i}\}_{i \in I}\,
      \rangle
    \end{align*}
    where:
    \begin{itemize}
    \item $\{\inj_{j}:Q_{j}\rTo Q\}_{j \in J}$ is a
      coproduct diagram in $\Cat{C}$,
    \item  $s_{i,j}: F_{i}Q_{j} \rTo
      Q$ for each $i \in I$ and $j \in J$,
    \item 
      $h_{i} : Q_{i} \rTo \Omega_{i}$ for each $i \in I$.
    \end{itemize}
  \item An arrow $f : \Aut{A} \rTo \Aut{B}$ is an arrow $f :
    Q^{\Aut{A}} \rTo Q^{\Aut{B}}$ in $\Cat{C}$ such that:
    \begin{itemize}
    \item for each $j \in J$, $\inj^{\Aut{A}}_{j}\comp f$
      factors (necessarily uniquely) through $\inj^{\Aut{B}}_{j}$:
      \begin{align*}
        \inj^{\Aut{A}}_{j}\comp f
        & = f_{j} \comp \inj^{\Aut{B}}_{j}\,,
      \end{align*}
    \item for each $i \in I$ and $j \in J$ the following
      equations hold:
      \begin{align*}
        s^{\Aut{A}}_{i,j}\comp f & = F_{i}f_{j}
        \comp s^{\Aut{B}}_{i,j}\,, &
        h^{\Aut{A}}_{i} & = f_{i} \comp h^{\Aut{B}}_{i} \,.
      \end{align*}
    \end{itemize}
  \end{itemize}
\end{definition}

\begin{lemma}
  \label{lemma:firstequivalence}
  The categories $\Pautomata$ and $\Pcoalg$ are equivalent.
\end{lemma}
\begin{proof}
  Given a $P$-coalgebra $\beta : B \rTo \sum_{i \in I} \Omega_{i}
  \times G_{i} B$, we construct a coproduct diagram $\inj_{i} : B_{i}
  \rTo B$ and arrows $\langle h_{i},t_{i}\rangle : B_{i} \rTo
  \Omega_{i} \times G_{i} B$ by pulling back along injections
  $\inj_{i} :\Omega_{i} \times G_{i} B \rTo \sum_{i \in I} \Omega_{i}
  \times G_{i} B$; we obtain a $P$-automaton $\Aut{B}$ by transposing
  the $t_{i}$.  This is the object part of a full and faithful functor
  $R :\Pcoalg \rTo \Pautomata$ that is the identity on morphisms.  In
  the other direction, given a $P$-automaton $\Aut{A}$ we define the
  coalgebra
  $$
  \sum_{i \in I} \langle h_{i},\be{s}_{i}\rangle : 
  Q^{\Aut{A}} \iso \sum_{i \in I} Q^{\Aut{A}}_{i}
  \rTo \sum_{i
    \in I} \Omega_{i} \times G_{i}Q^{\Aut{A}}\,.
  $$
  An arrow $f$ in $\Pautomata$ is also a coalgebra morphism, thus
  this construction defines a functor $L : \Pautomata \rTo \Pcoalg$
  which is left adjoint to $R$.
  
  $\Pcoalg$ is therefore a reflective subcategory of $\Pautomata$, and
  this depends only on distributivity of the base category $\Cat{C}$.
  In order to conclude that the two categories are equivalent it is
  must be argued that each object of $\Pautomata$ is isomorphic to an
  object coming from $\Pcoalg$; this is a consequence of coproduct
  injections being Cartesian natural transformations in an extensive
  category.
\end{proof}

In view of Lemma \ref{lemma:isos}, it should be easy to argue that the
category of $Q$-coalgebras is isomorphic to the category of
$F$-algebras equipped with an arrow from the carrier to $\Omega$, and
to the category of Eilenberg-Moore algebras for the monad $\F$
equipped again with an arrow from the carrier to $\Omega$. Using the
first isomorphism we are going to argue that:
\begin{lemma}
  \label{lemma:QAutequivQAlg}
  The category $Q$-coalgebras is equivalent to the category of
  $Q$-automata.
\end{lemma}  
\begin{proof}
  Given an $F$-coalgebra $s : FB \rTo B$ and an arrow $h : B \rTo
  \Omega$, we construct a coproduct diagram $\inj_{j} : B_{j} \rTo B$
  and arrows $h_{j} : B_{j} \rTo \Omega_{j}$ by pulling back along the
  coproduct diagram $\inj_{j} :\Omega_{j} \rTo \Omega$; we obtain a
  $Q$-automaton $\Aut{B}$ by letting $s_{i,j} = \inj_{F_{i}B_{j}}
  \comp F\inj_{j} \comp s$.  Again, this is the object part of a
  full and faithful functor that is the identity on morphisms. Its
  left adjoint is described as follows: given a $Q$-automaton
  $\Aut{A}$ we define $h = ! + \sum_{i \in I} h_{i} : Q^{\Aut{A}} \rTo
  \Omega$.  Recall that $(F_{i}\inj_{j}:F_{i}Q^{\Aut{A}}_{j} \rTo
  F_{i}Q^{\Aut{A}})_{j \in J}$ is a coproduct diagram, since $F_{i}$
  is a left adjoint. Therefore we define $s_{i}: F_{i}Q^{\Aut{A}} \rTo
  Q^{\Aut{A}}$ by saying that $F_{i}\inj_{j}\comp s_{i} = s_{i,j}$ and
  $s: FQ^{\Aut{A}} \rTo Q^{\Aut{A}}$ by $\inj_{i} \comp s = s_{i}$.
  
  The category of $F$-algebras equipped with an arrow to $\Omega$ is
  therefore a reflective subcategory of $\Qautomata$, assuming that
  $\Cat{C}$ is only distributive. The two categories are equivalent if
  $\Cat{C}$ is extensive as well.
\end{proof}

We are ready to define the functor $K$. In automata theoretic terms, it
amounts to completing a partial deterministic automaton by adding a
new sink state.
\begin{definition}
  For  a $P$-automaton $\Aut{A}$  the $Q$-automaton
  $K(\Aut{A})$ is defined as follows. We let $Q^{K(\Aut{A})}_{0}$ be the terminal
  object, and $Q^{K(\Aut{A})}_{i} = Q^{\Aut{A}}_{i}$ for $i\in I$; the
  coproduct diagram is given by:
  $$
  1 \rTo^{\inj_{l}} 1 + Q^{\Aut{A}} \lTo^{\inj_{r}} Q^{\Aut{A}}
  \lTo^{\inj^{\Aut{A}}_{i}}
  Q^{\Aut{A}}_{i} \,.
  $$
  The actions $s_{i,j}$ are defined as follows:
  \begin{align*}
    s_{i,j} & = \left \{
      \begin{array}{r@{\hspace{0mm}}ll}
        s^{\Aut{A}}_{i} \comp \inj_{r} &: F_{i}Q^{\Aut{A}}_{i} \rTo
        Q^{\Aut{A}} \rTo 1 + Q^{\Aut{A}}\,,
        & i = j\,,\\
        ! \comp \inj_{l} & :
        F_{i}Q^{\Aut{A}}_{j} \rTo 1 \rTo 1 + Q^{\Aut{A}} \,,
        &\textrm{otherwise}.
      \end{array}
    \right .
  \end{align*}
  Finally, the labeling is unchanged: $h^{K(\Aut{A})}_{i} =
  h^{\Aut{A}}_{i} $.  Clearly, if $f$ is a morphism in $\Pautomata$,
  then $1 + f$ is a morphism in $\Qautomata$, so that the construction
  $K$ defines a functor from $\Pautomata$ to $\Qautomata$.
\end{definition}

\begin{lemma}
  The functor $K$ is full and faithful.
\end{lemma}
\begin{proof}
  The functor is  faithful since $f$ is determined from $1 +
  f$ by pulling back along the right coproduct injection.  On the other
  hand, if $g : K(\Aut{A}) \rTo K(\Aut{B})$ is a morphism, then
  $\inj^{\Aut{A}}_{i} \comp \inj_{r} \comp g = g_{i} \comp
  \inj^{\Aut{B}}_{i} \comp \inj_{r}$, and $\inj_{l} \comp g = g_{0}
  \comp \inj_{l}$ where $g_{0}$ is necessarily the identity of the
  terminal object.  Therefore we can write $g = 1 + g'$, and $g'$ is
   a morphism in $\Qautomata$.  For example
   $$
   s^{\Aut{A}}_{i} \comp g' \comp \inj_{r}
   = 
   s^{\Aut{A}}_{i} \comp \inj_{r}\comp   g
   =
   s^{K(\Aut{A})}_{i,i} \comp   g
   =
   F_{i} g_{i} \comp s^{K(\Aut{B})}_{i,i} 
   =
   F_{i} g_{i} \comp s^{\Aut{B}}_{i} \comp \inj_{r}
   $$
  implies $s^{\Aut{A}}_{i} \comp g' = F _{i}g_{i} \comp
  s^{\Aut{B}}_{i}$, since coproduct injections are monic.
\end{proof}

\begin{definition}
  We denote by $\dQautomata$ the full subcategory category of
  $\Qautomata$ of those $Q$-automata $\Aut{B}$ for which
  $s^{\Aut{B}}_{i,i}$ factors through $Q^{\Aut{B}}_{I}$ and
  $s^{\Aut{B}}_{i,j}$ factors through $Q^{\Aut{B}}_{0}$ for $i \neq
  j$.
\end{definition}
The functor $K$ lands in $\dQautomata$, thus it can
be factored as
$$
\Pautomata \rTo^{K^{+}} \dQautomata \rTo^{I} \Qautomata\,,
$$
where the last functor $I$ is the inclusion. Since $K$ and $I$ are
both full and faithful, it follows that $K^{+}$ is full and faithful
too. Thus the image of the functor $K$ in the category $\Qautomata$
can be characterized as follows: an object $\Aut{B}$ of $\Qautomata$
is isomorphic to an object of the form $K(\Aut{A})$ if and only if
\begin{itemize}
\item it lies in the full subcategory $\dQautomata$,
\item $Q^{\Aut{B}}_{0}$ is a terminal object in $\Cat{C}$.
\end{itemize}
\begin{proposition}
  The functor $K^{+}$ has a left adjoint $L$.
\end{proposition}
The left adjoint $L(\Aut{A})$ is obtained by restricting all the
structure from $J$ to $I$: the underlying coproduct is now $\inj_{i}:
Q_{i} \rTo Q_{I}$; by the definition of $\dQautomata$ we can write
$s_{i,i} = s'_{i} \comp \inj_{I}$ (and such a factorization is unique)
so that the actions are defined to be these $s'_{i}$, while the
labeling arrows are unchanged.
%
Again, using the fact that coproduct injections are monic, it is
easily seen that this construction is functorial and that, for a
$Q$-automaton $\Aut{A}$ in $\dQautomata$ and a $P$-automaton
$\Aut{B}$, $f: Q^{\Aut{A}}_{I} \rTo Q^{\Aut{B}}$ is a morphism in
$\Pautomata$ if and only if $! + f: Q^{\Aut{A}} \rTo 1 + Q^{\Aut{B}}$
is a morphism in $\Qautomata$.



\begin{corollary}
  \label{cor:1}
  If the category $\dQautomata$ has finite limits, then the category
  $\Pautomata$ has also finite limits.
\end{corollary} 
The previous proposition has shown that we can identify $\Pautomata$
with a reflective subcategory of $\dQautomata$. It is therefore enough
to recall that a replete reflective subcategory is closed under
existing limits \cite[\S 3.5.3]{borceuxI}.

We shall argue in the next sections that $\dQautomata$ has indeed all
finite limits, so that the proposition holds without the proviso.


%% file: characterization.tex
\section{A Characterization of $\dQautomata$}
\label{sec:characterization}

In this section we propose an alternative characterization of
$Q$-automata in $\dQautomata$.  To this goal, we need a preliminary
observation: for an arbitrary object $\Aut{A}$ of $\Qautomata$ we
consider the two pullback squares:
$$
\mydiagram{
  []!W
  (!E{0.5}*+{\putequation{pb:argue}})
  !s{Q_{j}}{Q}{P_{j}}{\F Q}{1}{1}
  !a{^{\inj_{j}}}{^{\psi_{j}}}{^{\zero_{Q}}}{^{\pi_{\F Q}}}
  "3"
  !W
  (!E{0.5}*+{\putequation{pb:def}}) 
  !s{P_{j}}{\F
    Q}{Q_{j}}{Q}{1}{1} !a{^{\pi_{\F Q}}}{^{\pi_{j}}}{^{\s}}
  {|{\inj_{j}}} 
}\,.
$$
The square (\ref{pb:def}) is a pullback by definition. The square
(\ref{pb:argue}) is a pullback for the following reason. Because of
the unit law of Eilenberg-Moore algebras, the relation $\zero_{Q}
\comp \s = \id_{Q}$ holds. Also the commutative diagram corresponding
to the relation $\inj_{j} \comp \id_{Q} = \id_{Q_{j}} \comp \inj_{j}$
is a pullback since $\inj_{j}$ is monic. Hence (\ref{pb:argue}) is a
pullback and moreover $\psi_{j} \comp \pi_{j} = \id_{Q_{j}}$.  As a
consequence the diagram
$$
\mydiagram{
  []
  !s{F_{i}P_{j}}{F_{i}\F Q}
  {F_{i}Q_{j}}{Q}
  {1}{1}
  !a{^{F_{i}\pi_{\F Q}}}{^{F_{i}\pi_{j}}}
  {^{\mul_{i} \comp \s}}{^{s_{i,j}}}
}
$$
commutes, which can be seen as follows:
\begin{align*}
  F_{i}\pi_{\F Q} \comp \mul_{i} \comp \s 
  & =
  F_{i}\pi_{\F Q} \comp F_{i}\s \comp s_{i} 
  = F_{i}(\pi_{\F Q} \comp \s) \comp
  s_{i}  \\
  & = \by{pb:def} 
  F_{i}(\pi_{j} \comp \inj_{j}) \comp s_{i}  \\
  & = 
  F_{i}\pi_{j}  \comp s_{i,j} \,. 
  \tag*{\text{by Lemma \ref{lemma:QAutequivQAlg}}} 
\end{align*}
Considering that for an object in $\dQautomata$ we have $s_{i,i} =
s'_{i,i} \comp \inj_{I}$ and $s_{i,j} = s'_{i,j} \comp \inj_{0}$ if $i
\neq j$, we see that the diagrams
$$
\mydiagram{
  [](!E{0.8}
  []*+{F_{i}P_{i}}
  (:[r(0.8)]*+{F_{i}Q_{i}}^{F_{i}\pi_{i}}
  :[r(0.8)]*+{Q_{I}}^{s'_{i,i}}
  :[d]*+{Q}="G"^{\inj_{I}})
  :[d]*+{{F_{i}\F Q }}|{ F_{i}\pi_{\F Q}  }
  :"G"^{\mul_{i} \comp \s}
}
\quad\quad
\mydiagram{
  [](!E{0.8}
  []*+{F_{i}P_{j}}
  (:[r(0.8)]*+{F_{i}Q_{j}}^{F_{i}\pi_{j}}
  :[r(0.8)]*+{Q_{0}}^{s'_{i,j}}
  :[d]*+{Q}="G"^{\inj_{0}})
  :[d]*+{{F_{i}\F Q }}|{ F_{i}\pi_{\F Q}  }
  :"G"^{\mul_{i} \comp \s}
}
\;\;
i\neq j
$$
commute.  If we define
\begin{align*}
  P_{i,j} &
  = \text{pullback of } \mul_{i} \comp \s  
  \text{ against }
  \begin{cases}
    \inj_{I} & \text{if } i = j\,, \\
    \inj_{0} & \text{otherwise}\,,
  \end{cases}
\end{align*}
we obtain that the arrow $\F_{i} \pi_{\F Q}$ factors through the
subobject $P_{i,j}$ of $F_{i} \F Q$. Transposing this relation
according to Lemma \ref{lemma:univquantification}, we deduce that, for
each $j \in J$ and $i \in I$, the relation
\begin{align}
  \label{rel:internal}
  P_{j} & \leq \forall_{F_{i}} P_{i,j}
\end{align}
holds in the Brouwerian semilattice of subobjects of $\F Q$.
Conversely, suppose that $P_{j} \leq \forall_{F_{i}} P_{i,j}$, that is
we can write $F_{i} \pi_{j} \comp s_{i,j} = F_{i}\pi_{\F Q} \comp
\mul_{i} \comp \s = \alpha \comp \inj_{d}$, where $d = I$ if $i = j$
and $d = 0$ otherwise. It follows that $s_{i,j} = F_{i}\psi_{j} \comp
F_{i} \pi_{j} \comp s_{i,j} = F_{i}\psi_{j} \comp \alpha \comp
\inj_{d}$, which shows that $s_{i,j}$ factors through the proper
coproduct injection.  Thus we have shown that:
\begin{proposition}
  \label{prop:characterization}
  An object $\Aut{A}$ of $\Qautomata$ belongs to $\dQautomata$ if and
  only if, for each $i \in I$ and $j \in J$, the relation
  \eqref{rel:internal} holds in the Brouwerian semilattice of
  subobjects of $\F Q^{\Aut{A}}$.
\end{proposition}


%% file: specialQaut.tex
\section{The Coreflector}
\label{sec:coreflection}

The observations of the previous section suggest the following
construction to be performed on an arbitrary $Q$-automaton $\Aut{A}$.
Its building blocks are subobjects $P_{i}$ of $\F Q^{\Aut{A}}$ and
$P_{i,j}$ of $F_{i}\F Q^{\Aut{A}}$ defined as the following pullbacks:
$$
\mydiagram[6em]{
  []!W
  !s{P_{j} }
  {\F Q  }{Q_{j}}{Q}
  {1}{1.2}
  !a{}{}{^{\s}}{^{\inj_{j}}}
}
$$
$$
\mydiagram[6em]{
  [](
  []
  )
  !W
  !s{ P_{i,i} }
  {F_{i}\F Q  }{Q_{I}}
  {Q}
  {1}{1.2}
  !a{}{}{^{\mul_{i} \comp \s}}{^{\inj_{I}}}
}
\quad
\mydiagram[6em]{
  [](
  []
  )
  !W
  !s{ P_{i,j} }
  {F_{i}\F Q}{Q_{0}}{Q}
  {1}{1.2}
  !a{}{}{^{\mul_{i} \comp \s}}{^{\inj_{0}}}
}
\;\;i \neq j \,.
$$
The $P_{j}$ and the $P_{i,j}$ are indeed subobjects of $\F Q$ and
$F_{i}\F Q $, respectively, because of Lemma \ref{lemma:robin}.  In
the Brouwerian semilattice of subobjetcs of $Q^{\Aut{A}}$, we define:
\begin{align*}
  C_{j} 
  & =
  \forall_{\F}(\,
  P_{j} \rightarrow \bigwedge_{i \in I} \forall_{F_{i}} P_{i,j}
  \,)\,,
  &
  D 
  & = \bigwedge_{j \in J} C_{j}\,.
\end{align*}
The meaning of the universal quantification $\forall_{F}$ and the
reason for which $C_{j}$ is a subobject of $Q^{\Aut{A}}$ are explained
in Lemma \ref{lemma:univquantification}.

In order to understand the definition of the $C_{j}$, we consider the
set theoretic Example \ref{ex:main}. Here $Q = \Omega^{A^{\ast}}$ is
the set of all $A$-complete trees labeled either by the symbol $\bot$
or by some element in some $\Omega_{i}$; $\s :A^{\ast} \times
\Omega^{A^{\ast}} \rTo \Omega^{A^{\ast}}$ is correspondence which
takes the pair $(w,t)$ to the subtree of $t$ rooted at $w$, the
function $\lambda x.t(wx)$, and $h : \Omega^{A^{\ast}} \rTo \Omega$ is
evaluation at the empty word.  $P_{i,i}$ is seen to be the set of
triples $(a,w,t)$ in $A_{i}\times A^{\ast} \times \Omega^{A^{\ast}}$
such that $t(wa)$ belongs to $\Omega_{i}$ for some $i \in I$, and
$P_{i,j}$ -- for $i \neq j$ -- is the set of triples $(a,w,t)$ in
$A_{j}\times A^{\ast} \times \Omega^{A^{\ast}}$ such that the subtree
of $t$ rooted at $wa$ is labeled by $\bot$.  Therefore $C_{j}$ is the
collection of trees $t \in \Omega^{A^{\ast}}$ with the following
property: \emph{for all $w \in A^{\ast}$ such that $t(w) \in
  \Omega_{j}$, for all $i \in I$ and $a \in A_{i}$, if $i = j$ then
  $t(wa)$ is in $\Omega_{I}$, and otherwise, if $i \neq j$, then
  $t(wa) = \bot$.}  That is, the definition of $C_{j}$ mimics clauses
(\ref{clause:f})-(\ref{clause:g}) of \ref{ex:main}.

\begin{lemma}
  \label{lemma:first}
  The object $D$ is a subautomaton of $\Aut{A}$. 
\end{lemma}
\begin{proof}
  We shall show that the arrow $\F \inje{C_{j}} \comp \s$ factors
  through $C_{j}$; this will be enough since it is easily verified
  that the intersection of subobjects closed under the action of $\F$
  is again closed under this action. By the definition of $C_{j}$,
  this is equivalent to the pullback of $P_{j}$ along $\F(\F
  \inje{C_{j}} \comp \s)$ to factor through $\forall_{F_{i}} P_{i,j}$,
  for all $i \in I$.  To this end, in the diagram below the arrow
  $\psi$ corresponds under several adjunctions to the identity of
  $C_{j}$. This implies that we can factor the pullback $P$ through
  $\forall_{F_{i}} P_{i,j}$:
  $$
  \mydiagram{
    []!W
    []!s{P}{\F \F C_{j}}{\F C_{j} \pb[1.5]{\F Q  } P_{j}}
    {\F C_{j}}
    {1}{1.5}
    !a{}{}{^{\mul_{C_{j}}}}{}
    "3"!W[]"3"
    !s{\F C_{j} \pb[1.5]{\F Q  } P_{j}}{\F C_{j}}{P_{j}}{\F Q  }
    {1}{1.5}
    !a{}{}{^{\F \inje{C_{j}}}}{}
    "1"="1m"
    "4"="FT"
    "1m":@/^2em/[r(2.2)]*+{\forall_{F_{j}} P_{i,j}}^{\psi}:"FT"
  }
  $$
  The statement of the lemma  follows since  the relations
  $$
  \mul_{C_{j}} \comp \F \inje{C_{j}} \comp \s
  = \F \F \inje{C_{j}} \comp \mul_{Q  } \comp  \s \\
  = \F \F \inje{C_{j}} \comp \F \s \comp \s
  $$
  exhibit $P$ as the pullback of $P_{j}$ along $\F(\F \inje{C_{j}}
  \comp \s)$.
\end{proof}

\begin{proposition}
  The automaton $D$ belongs to the category $\dQautomata$.
\end{proposition} 
\begin{proof}
  The diagram
  $$
  \mydiagram{ 
    []="S"!Z 
    [dr] !W !s{\F C_{j} \pb[1.5]{\F Q }
      P_{j}}{\F C_{j}}{P_{j}}{\F Q } {1}{1} 
    !a{}{}{^{\F\inje{C_{j}}}}{}
    "1"="1m" "2"="2m"[l]="2p" 
    "3" !W !s{P_{j}}{\F
      Q }{Q_{j}}{Q} {1}{1} 
    !a{}{}{^{\s}}{^{\inj_{j}}} "2"="FT"
    "3"
    [u]="3p" 
    "S"*+{D_{j}} (:"3p"*+{Q_{j}}:@2@{-}"3",
    :"2p"*+{D}="2p"|{\inj_{j} = \inje[D]{D_{j}}} :"2m"^{\zero_{D} \comp
      \F\inje[C_{j}]{D}}, :"1m" ) 
    "2p":[r(0.5)u(-0.5)]*+{Q
    }_{\inje[Q]{D}}:"FT"|{\zero_{Q }}
    "1m":@/^2em/[r(1.5)u(0.2)]*+{\forall_{F_{i}} P_{i,j}}^{\psi} :"FT"
    "Q ":"4"_{\id_{Q}} }
  $$
  shows that the arrow $\inje[Q]{D_{j}} \comp \zero_{Q }$ can be
  factored through $\forall_{F_{i}}P_{i,j}$ over $\F Q $. Transposing
  this relation, we obtain that $F_{i}(\inje[Q]{D_{j}} \comp \zero_{Q
  })$ can be factored through $P_{i,j}$ over $F_{i}\F Q $.
  Considering the definition of the $P_{i,j}$ as pullbacks, we obtain
  that $F_{i}(\inje[Q]{D_{j}} \comp \zero_{Q})\comp \mul_{i} \comp \s$
  can be factored through $\inj_{I}: Q_{I} \rTo Q$ if $i = j$ and
  through $\inj_{0}: Q_{0} \rTo Q$ otherwise.
  Observe that $F_{i}\zero_{Q }\comp \mul_{i} \comp \s = s_{i}$ and
  that $F_{i}(\inje[Q]{D_{j}}\comp \zero_{Q }) \comp \mul_{i} \comp \s
  = F_{i}\inje[Q]{D_{j}} \comp s_{i} = s^{D(\Aut{A})}_{i,j} \comp
  \inje[Q]{D}$, thus: $s^{D(\Aut{A})}_{i,j} \comp \inje[Q]{D}$ can be
  factored through $\inj_{I}: Q_{I} \rTo Q$ if $i = j$ and through
  $\inj_{0}: Q_{0} \rTo Q$ otherwise.  Since $D_{0}$ and $D_{I}$ are
  obtained by pulling back $Q_{0}$ and $Q_{I}$ against $\inje[Q]{D}$,
  we obtain the statement of the proposition.
\end{proof}

\begin{proposition}
  Let $\Aut{A},\Aut{B}$ be two automata, of which $\Aut{A}$ is in
  $\dQautomata$. If $f : \Aut{A} \rTo \Aut{B}$ is a morphism, then it
  factors through $D(\Aut{B})$.
\end{proposition}
\begin{proof}
  It is an easy exercise to show that if $\inje{} : \Aut{C}
  \rTo \Aut{B}$ is monic, and the morphism $f: \Aut{A} \rTo \Aut{B}$
  factors through $Q^{\Aut{C}}$ in the underlying category $\Cat{C}$,
  $f = f' \comp \inje{}$, then $f'$ is also a morphism from $\Aut{A}$
  to $\Aut{C}$ in $\Qautomata$.  
  
  Therefore we shall prove that $f$ factors through $D =
  Q^{D(\Aut{B})}$, by showing that it factors through each
  $C^{\Aut{B}}_{j}$.  Unraveling its definition, we need to show that
  $\F f$ factors through $P^{\Aut{B}}_{j} \rightarrow \bigwedge_{i\in
    I}\forall_{F_{i}} P^{\Aut{B}}_{i,j}$, or equivalently that the
  pullback $P_{j}$ of $\F f \comp \s^{\Aut{B}}$ along the injection
  $Q^{\Aut{B}}_{j} \rTo Q^{\Aut{B}}$ factors through each
  $\forall_{F_{i}} P^{\Aut{B}}_{i,j}$ over $\F Q^{\Aut{B}}$. Since $\F
  f \comp \s^{\Aut{B}} = \s^{\Aut{A}} \comp f$ and $h^{\Aut{A}} = f
  \comp h^{\Aut{B}}$, this pullback is $P^{\Aut{A}}_{j}$.  Since
  $\Aut{A}$ belongs to $\dQautomata$, we know from Proposition
  \ref{prop:characterization} that $P^{\Aut{A}}_{j}$ factors through
  $\forall_{F_{i}} P^{\Aut{A}}_{i,j}$ over $\F Q^{\Aut{A}}$.  In order
  to reach our goal, we only need to argue that the construction
  $\forall_{F_{i}} P^{\Aut{X}}_{i,j}$ is natural in $\Aut{X}$.
  Clearly the $P^{\Aut{X}}_{i,j}$ are natural in $\Aut{X}$, and
  transposing the diagram
  $$
  \mydiagram{
    []="S"
    [r]!s{P^{\Aut{A}}_{i,j}}{F_{i}\F Q^{\Aut{A}}}
    {P^{\Aut{B}}_{i,j}}{F_{i}\F Q^{\Aut{B}}}
    {1}{1}
    !a{}{}{^{F_{i}\F f}}{}
    "S"
    *+{F_{i}\forall_{F_{i}}P^{\Aut{A}}_{i,j}}
    (:"1",:"2")
  }
  $$
  whose triangle on the left is the counit of the adjunction, we
  derive the factorization we are looking for.
\end{proof}

To conclude this section, we resume what we have proved:
\begin{theorem}
  The category $\dQautomata$ is a coreflective subcategory of
  $\Qautomata$, the coreflector being the construction $D$.
\end{theorem}

\begin{corollary}
  \label{cor:2}
  If $\Qautomata$ has finite limits, then $\dQautomata$ has finite
  limits.
\end{corollary}
It is a standard fact \cite[\S 3.5.3]{borceuxI} that the limit of a
diagram $\set{\Aut{A}_{k}}_{k \in K}$ in $\dQautomata$ can be
calculated by applying the coreflection to the $Q$-automaton $\lim_{k}
\Aut{A}_{k}$.

We can now argue that $\Qautomata$ has all finite limits, so that the
statements of Corollaries \ref{cor:1} and \ref{cor:2} hold without
their proviso.
\begin{proposition}
  The category $\Qautomata$ has all finite limits.
\end{proposition}
\begin{proof}
  We have argued that $\Qautomata$ is isomorphic to the category of
  Eilenberg-Moore algebras for the free monad on $F$ equipped with an
  arrow to $\Omega$. Since we are assuming that an $F$-algebra, call
  it $\G \Omega$, cofree over $\Omega$ exists, we end-up 
  establishing the following isomorphism of categories:
  \begin{align*}
    \Qautomata & \iso \Cat{C}^{\bbt}/\G \Omega\,,
  \end{align*}
  where $\bbt = (\F,\zero,\mul)$ is the free monad on $F$.  That is,
  we have identified $\Qautomata$ with a certain slice category.
  Existence of finite limits follows then from existence of finite
  limits in $\Cat{C}^{\bbt}$ which are created from the forgetful
  functor.
\end{proof}


%% file: further.tex
\section{Further Observations}
\label{sec:further} 

\subsection{Behavioral Characterization of $\dQautomata$.}

The category $\dQautomata$ has the following property: if $\Aut{B}$ is
in $\dQautomata$ and $f: \Aut{A}\rTo \Aut{B}$ is an arrow of
$\Qautomata$, then $\Aut{A}$ is also in $\dQautomata$. This can be
seen by considering the following diagram showing that if
$s^{\Aut{B}}_{i,i}$ factors through $\inj_{I}$, then so does
$s^{\Aut{A}}_{i,i}$ (a similar diagram is used for $s_{i,j}$ when $i
\neq j$):
$$
\mydiagram{
  []*+{F_{i} Q^{\Aut{A}}_{i}}="Start"
  [r(1)d(0.6)]!W
  !s{Q_{I}^{\Aut{A}}}{Q^{\Aut{A}}}
  {Q_{I}^{\Aut{B}}}{Q^{\Aut{A}}}
  {1}{1}
  !a{|{\inj_{I}}}{^{f_{I}}}{^{f}}{|{\inj_{I}}}
  "Start"(
  :[r(2)]*+{F_{i} Q^{\Aut{B}}_{i}}="UR"^{F_{i}f_{i}}
  :"3"|{s_{i,i}'},
  :@."1",
  :"2"^{s_{i,i}^{\Aut{A}}}
  )
  "UR":@/^7mm/"4"^{s_{i,i}^{\Aut{B}}}
}
$$
Hence we can identify $\dQautomata$ with the slice category
$\Qautomata /D(\G \Omega)$:
\begin{proposition}
  A $Q$-automaton belongs to $\dQautomata$ if and only if the unique
  arrow to the terminal $Q$-automaton $\G \Omega$ factors through
  $D(\G \Omega)$.
\end{proposition}
\noindent
In particular, the coreflection $D$ can be described as pulling back
along the monomorphism $\inje[\G \Omega]{D(\G \Omega)}: D(\G \Omega)
\rTo \G \Omega$. If we think of the unique arrow from a $Q$-automaton
to the terminal $Q$-automaton as its behavior, we see that the
category $\dQautomata$ is completely determined by the generic
behavior $D(\G \Omega)$.  This discussion also shows that describing
the category $\dQautomata$ as a coreflective subcategory of
$\Qautomata$ and describing the object $D(\G \Omega)$ are under some
extent equivalent. The latter approach was the one pursued in an
earlier version of this paper \cite{previous}.

\subsection{Extensiveness is Necessary}

Finally, the reader should notice that all over section
\ref{sec:coreflection} we have implicitly used extensiveness of
$\Cat{C}$: we have assumed that every $Q$-automaton comes with an
associated $F$-algebra labeled by $\Omega$ so that diagrams such as
$$
\mydiagram{
[]!s{Q_{i}^{\Aut{A}}}{Q^{\Aut{A}}}{\Omega_{i}}{\Omega}
{1}{1}
!a{^{\inj_{i}}}{^{h_{i}}}{^{h}}{^{\inj_{i}}}
}
\quad\quad
\mydiagram{
[]!s{Q_{j}^{\Aut{A}}}{Q^{\Aut{A}}}{Q_{j}^{\Aut{B}}}{Q^{\Aut{B}}}
{1}{1}
!a{^{\inj_{j}}}{^{f_{j}}}{^{f}}{^{\inj_{j}}}
}
$$
are pullbacks. 
  If extensiveness were not used at this point, we could use the diagram of reflective and coreflective subcategories
$$
\newdir{ (}{\dir{}!/-4pt/\dir^{(}}
\newdir{ )}{\dir{}!/-7pt/\dir_{(}}
\mydiagram[7em]{
  []*+{\Pcoalg}="1"
  :@{{ (}->}[r]*+{\Pautomata}="2"^{}="A"
  (:@/_1.3em/"1"^{}="B")
  ("A":@{{}{}{}}"B"|{\bot})
  :@{ (->}[r]*+{\dQautomata}="3"^{}="A"
  (:@/_1.3em/"2"^{}="B")
  ("A":@{{}{}{}}"B"|{\bot})
  :@{ (->}[r]*+{\Qautomata}="4"^{}="A"
  (:@/^1.3em/"3"^{}="B")
  ("A":@{{}{}{}}"B"|{\bot})
  :@/^1.3em/[r]*+{\Qcoalg}^{}="A"
  (:@{ )->}"4"^{}="B")
  ("A":@{{}{}{}}"B"|{\bot})
}
$$
to argue that the existence of a terminal object in $\Qcoalg$ gives
rise to a terminal object in $\Pcoalg$. It is not the case, however,
that the construction $D(\G \Omega)$ gives rise to a final
$P$-coalgebra if coproducts are not disjoint. In this sense,
extensiveness of $\Cat{C}$ is a necessary condition. For example, if
we consider an arbitrary Heyting algebra, we need to ask whether the
relation
\begin{align}
  \label{eq:segerberg}
  \nu_{X}.(\, \bigvee_{i \in I} \,(\,\Omega_{i} \land G_{i}X\,)\, ) &
  = \Omega_{I} \land \bigwedge_{i \in I} \G(\Omega_{i}\rightarrow
  G_{i}\Omega_{I})
\end{align}
holds, that is, whether the expression on the right -- which is the
posetal version of the pullback of $\Omega_{I}$ against $D(\G \Omega)$
-- produces a greatest fixed point of the function $\bigvee_{i}
\Omega_{i} \land G_{i}X$.  When the indexing set $I$ is a singleton
this equation holds.\footnote{It amounts to Segerberg's axiomatization
  of Propositional Dynamic Logic \cite{segerberg}} On the other hand,
the transition system below provides a counterexample to the binary
version:
$$
\mydiagram[3em]{
  []*+{a}="A"
  :[r(2.5)]*+{b}="B"^{1}
  "A"[u(0.5)]*+{
  }
  "B"[u(0.5)]*+{\Omega_{1},\Omega_{2}}
}
$$
Here, the transition system has two states, $a$ and $b$, and a
transition from $a$ to $b$  labeled by $1$.  Over this transition
system we interpret the modal operator $\pos[1]$, by means of the
relation labeled by $1$, and a second modal operator $\pos[2]$ by
means of the empty relation. We also declare that the two
propositional constants $\Omega_{1}$ and $\Omega_{2}$ hold in state
$b$ only.  If we put $\pos[\set{1,2}^{\ast}]Y = \mu_{X}.(Y \vee
\pos[1]X \vee \pos[2]X)$, then we observe that the relation
\begin{align*}
  \makebox[4cm][l]{$(\Omega_{1} \land \Omega_{2}) \;\vee\;
    \pos[\set{1,2}^{\ast}]( \,\neg
    \Omega_{1}\, \land \,\pos[1](\Omega_{1} \land \Omega_{2})\, ) $} \\
  \makebox[1cm][l]{$\vee\; \pos[\set{1,2}^{\ast}] ( \,\neg \Omega_{2}
    \,\land \,\pos[2](\Omega_{1} \land \Omega_{2})\, )$}
  \\
  &\quad \quad\leq \mu_{X}.(\, (\Omega_{1} \vee \pos[1]X) \,\land\,
  (\Omega_{2} \vee \pos[2]X) \,)
\end{align*}
does not hold in the transition system: this relation is the dual of
(\ref{eq:segerberg}).


%% file: ack.tex
\section*{Acknoweledgment}
The author is grateful to professor Robin Cockett for useful
discussion on the subject.
